\documentclass[]{interact}
\usepackage[T1]{fontenc} 
\usepackage[utf8]{inputenc} 
\usepackage{floatrow} 
\makeatletter
\renewcommand*\l@section{\@dottedtocline{1}{1.5em}{2.3em}}
\makeatother\usepackage{tablefootnote}
\usepackage[linktoc=page]{hyperref}
\hypersetup{linkcolor  = black}
\hypersetup{citecolor=black} 
\hypersetup{colorlinks=true}
\hypersetup{urlcolor=cyan}
\usepackage{color, colortbl}
\usepackage{graphicx, amsmath, amsthm, amssymb, mathrsfs, amscd}
\usepackage{tablefootnote}
\usepackage{footnote}
\usepackage{epstopdf}
\usepackage[caption=false]{subfig}
\usepackage[numbers,sort&compress]{natbib}
\bibpunct[, ]{[}{]}{,}{n}{,}{,}
\makeatletter
\def\NAT@def@citea{\def\@citea{\NAT@separator}}
\makeatother
\theoremstyle{plain}
\newtheorem{theorem}{Theorem}[section]
\newtheorem{lemma}[theorem]{Lemma}

\theoremstyle{definition}
\newtheorem{definition}[theorem]{Definition}

\theoremstyle{remark}

\usepackage{lineno,hyperref}
\modulolinenumbers[5]
\usepackage{tablefootnote}
\hypersetup{linkcolor  = blue}
\hypersetup{citecolor=blue} 
\hypersetup{colorlinks=true}
\hypersetup{urlcolor=cyan}
\usepackage{color, colortbl}
\usepackage{multirow}
\usepackage{graphicx, amsmath, amsthm, amssymb, mathrsfs, amscd}
\usepackage{caption}
\usepackage{enumerate}
\usepackage{algorithm}
\usepackage[noend]{algcompatible}
\usepackage{url}
\usepackage{tikz}
\usepackage[toc,page]{appendix}
\usepackage{footnote}

\usepackage{color}
\usepackage{color, colortbl}
\usepackage{tablefootnote}
\usepackage{footnote}
\theoremstyle{definition}

\theoremstyle{definition}

\theoremstyle{definition}

\theoremstyle{definition}
\theoremstyle{definition}
\makeatletter
\def\cleardoublepage{\clearpage\if@twoside \ifodd\c@page\else
  \hbox{}
  \vspace*{\fill}
    \vspace{\fill}
  \if@twocolumn\hbox{}\newpage\fi\fi\fi}
\makeatother

\begin{document}
\title{Curves formed by Vanishing Discriminant and Roots of Complex-valued Harmonic Polynomials\\ ~~~~~~~~~~~~~~~~~~~(Computer-Aided  Case Study)}
\author{
\bigskip \name{Oluma~Ararso~Alemu \textsuperscript{1,\textcolor{blue}{*}} and Hunduma~Legesse~Geleta \textsuperscript{1}}\footnote{ \textcolor{blue}{*}Corresponding author}
 \affil{\textsuperscript{1} Department of Mathematics, Addis Ababa University, Addis Ababa, Ethiopia.}}
\maketitle
\begin{abstract}
In this paper, we determine and specify the type of curves formed by the vanishing discriminant of some specified family of complex-valued harmonic polynomials with two parameters. We also classify the region formed by curves as bounded and unbounded connected components which in turn used to count the zeros of the complex-valued harmonic polynomials. Our study is a computer-aided case and the result shows that the curves formed are a teacup curve, a parabolic curve and a  swallowtail catastrophe curve. These curves come up together to form a butterfly catastrophe.
\end{abstract}
\section{Introduction}

The study of the zero of analytic polynomials  has become of interest due to their application in other fields. The Fundamental Theorem of Algebra asserts  that a polynomial of degree $n\geq 1$ has at most  $n$ real roots and exactly $n$ complex roots counting with multiplicities. The Fundamental Theorem of Algebra holds only in univariate and analytic polynomials. If a polynomial is multivariate or non-analytic Fundamental Theorem of Algebra does not hold in general. In the plane, Bezout's Theorem states that the number of common zeros of system of polynomials equals the  product of the \bigskip degrees of polynomials. \\
In studying zeros of complex-valued harmonic polynomials, we found the zero inclusion regions of general complex-valued harmonic polynomials in  \cite{legesse2022location}. In the same paper, we solved an interesting problem stated in \cite{brilleslyper2020zeros} by Brilleslyper et al. We found the zero inclusion region of complex-valued harmonic trinomials.  Continuing our work we  considered, in \cite{alemu2022zeros} and  \cite{ararso2022image}, the harmonic quadrinomial of the form $Q_{b,c}(z)= bz^k+\overline{z}^n+c\overline{z}^m+z$  with $b,c \in \mathbb{R}\setminus\{0\},$ $k\geq n>m,$ and $k,n,m \in \mathbb{N}$. As can be seen, paper \cite{alemu2022zeros} is an updated version of the preprint \cite{ararso2021zeros} and we have determined the critical curve separating sense-preserving region from sense-reversing region. The main motivation for the present paper comes from the study of the zeros of complex-valued harmonic polynomials. The main interest here is geometrical information about the curves formed by the vanishing discriminant of systems of polynomials formed from the complex-valued harmonic polynomials  and counting the number of zeros in each connected components \bigskip formed by the curves of the discriminant.\\ 
It is well known that a catastrophe theory is a branch of bifurcation theory in the study of dynamical systems; it is also a particular special case of more general singularity theory in geometry. Our study in this paper is related to a catastrophe theory, since it also studies how the qualitative nature of equation solutions depends on the parameters that appear in the equations. Applying a vanishing discriminant to quadrinomials $\overline{z}^4=z^2+az+b,$ we find that the curves formed are parabolic curve, teacup curve and swallowtail catastrophe curve. When studying these  families of harmonic polynomials in $\mathbb{C}$ and looking at the curves of the vanishing discriminant using Mathematica, the \bigskip curve formed is also butterfly catastrophe.\\
This paper is organized as follows. In section \ref{p}, we present some important preliminary results that will formalizes the main results. In section $\ref{c},$ we determine and specify the types of curves formed by the vanishing discriminant of the given harmonic polynomial and we classify the curve formed as bounded and unbounded connected components. In section $\ref{e},$ we count the zeros of $f(z)= \overline{z}^4-(z^2+az+b)$ for $a,b \in \mathbb{R}$ in both bounded and unbounded components of the curves formed by the vanishing discriminant. As a result,  $4 \leq \mathscr{Z}_f \leq 6$ for  unbounded connected components and $4 \leq \mathscr{Z}_f \leq 10$ for bounded connected components, where $\mathscr{Z}_f$ denotes the number of zeros of $f.$ In section $\ref{r}$ we have some important remarks and conclusions.
\section{Preliminaries} $\label{p}$
In this section we review some important concepts and results that we will use later in the main result. We begin by stating the well known results and some useful definitions, theorems and lemmas. More specifically, we focus on discriminant, resultant, some results on both and theorems on bounding the zeros.\\

The discriminant of a polynomial is a polynomial function of the coefficients of the original polynomial. It is a quantity that depends on the coefficients and allows deducing some properties of the roots without computing them. We usually use discriminant in polynomial factoring, number theory and algebraic geometry. The discriminant of a polynomial gives some insight into the nature of the zeros of a polynomial.\\

The discriminant of a quadratic polynomial  $p(x)=ax^2+bx+c$ is given by 
\begin{equation}
\mathcal{D}_2 = b^2-4ac.
\end{equation} Here, if $\mathcal{D}_2 >0,$ then $p(x)$ has two distinct real roots. If $\mathcal{D}_2 =0,$ then $p(x)$ has one real root with multiplicity two. If $\mathcal{D}_2<0,$ then $p(x)$ has no real root, but it has two complex conjugates roots. It has been known since the sixteenth century that a cubic polynomial $p(x)=ax^3 + bx^2 + cx + d$ has a repeated root if and only if its discriminant, 
\begin{equation}
\mathcal{D}_3=b^2c^2-4ac^3-4b^3d-27a^2d^2+18abcd
\end{equation} is zero. Continuing like this, the discriminant of polynomial of any degree $n\geq 1$ is defined as follows in general.

\begin{definition}\cite{janson2007resultant}
The discriminant of the polynomial, $p(z)=a_nz^n+a_{n-1}z^{n-1}+ \cdots +a_2z^2+a_1z+a_0,$ of degree $n$ with roots $r_1, r_2, \cdots r_n$ is defined as 
\begin{equation}
\mathcal{D} _n = a_n^{2(n-1)}\prod_{i,j~~i<j}^n(r_i-r_j)^2,
\end{equation}
which gives a homogeneous polynomial of degree $2(n-1)$ in the coefficients of $p.$
\end{definition} 
Recall that the discriminant of a univariate polynomial of positive degree is zero if and only if the polynomial has a multiple roots.  For a polynomial with real coefficients with no multiple roots, the discriminant is positive if the number of non-real roots is \bigskip a multiple of 4 and negative  otherwise.\\
To determine the existence of a root to  a system of polynomial equations, the resultants are applicable and are used to reduce a given system to one with fewer variable. The resultant of two univariate polynomials is also used to decide whether they have common zero(this works efficiently for any polynomials).
\begin{definition}\cite{woody2016polynomial}
Given two polynomials $f(z)=a_0+a_1z+a_2z^2+a_z^3+ \cdots + a_nz^n $ and $g(z)=b_0+b_1z+b_2+ \cdots + b_mz^m$ over $\mathbb{C}.$ The resultant of $f$ and $g$ relative to the variable $z$ is a polynomial over the field of coefficients of $f(z)$ and $g(z);$ and is defined as 
\begin{equation}
 Res(f,g,z)=a_n^mb_m^n\prod_{i,j}(\alpha _i-\beta _j) 
 \end{equation} where $f(\alpha _i)=0$ for all $1 \leq i \leq n$ and $g(\beta _j)=0$ for all $1 \leq j \leq m.$
\end{definition}
The following lemma is proved in \cite{woody2016polynomial}.
\begin{lemma}
The resultant of $f(x)$ and $g(x)$ is equal to zero if and only if the two polynomials have a root in common.
\end{lemma}
The vanishing discriminant of a polynomial $f$ and the resultant of $f$ with $f'$ have a nice relationship and is illustrated in the following lemma.
\begin{lemma}
The discriminant of a polynomial $f$ over its domain vanishes if and only if $Res(f,f',z)=0.$
\end{lemma}
In any simply connected sub-domain of $\mathscr{G}\subset \mathbb{C}$  a complex-valued harmonic polynomial $f(z)=u(x,y)+iv(x,y)$ can be decomposed as $ f(z) = h(z) + \overline{g(z)},$  where both $g$ and $h$ are analytic polynomials. In this decomposition,  $h(z)$ is called analytic  part and $g(z)$ is said to be co-analytic part of $f.$ This family of complex-valued harmonic functions is a generalization of analytic mappings studied in geometric  function theory, and much research has been  done  investigating the properties of these harmonic \bigskip functions.\\
Wilmshurst \cite{hauenstein2015experiments} considered such complex-valued harmonic polynomials of the form  
$f(z)=h(z)+\overline{g(z)}.$ If $\mathrm{deg}h=n>m=\mathrm{deg}g,$ then as to the question of improving the bound $\mathscr{Z}_f \leq n^2$ given additional information, Wilmshurst made the conjecture $\mathscr{Z}_f \leq 3n-2+m(m-1).$ This conjecture is stated in \cite{wilmshurst1998valence}. It is also among the list of open problems in \cite{bshouty2010problems}. For $m=n-1$ the upper bound follows from Wilmshurst's theorem  and examples were also given in  \cite{wilmshurst1998valence} showing that this bound is sharp. For $m = 1,$ the upper bound was proved by D.Khavinson and   G.Swiatek \cite{khavinson2003number}, and bound was also sharp.  For $m = n- 3$, the conjectured bound is $3n - 2+ m(m- 1) = n^ 2 - 4n + 10.$ It was lso shown that  $\mathscr{Z}_f > n^ 2 - 3n + \mathcal{O}(1).$ For this lower bound counterexample was \bigskip given in \cite{lee2015remarks}.

\begin{definition}
The valence of a function $f$ at a given point $\omega,$ denoted by $Val(f, \omega),$ is the number of distinct points $z$ in the domain of $f$ such that $f(z) = \omega.$ The valence of a function, denoted by $Val(f),$ is the supremum of $Val(f, \omega),$ for each $\omega$ in the domain of $f.$
\end{definition}
\begin{definition}
  A harmonic function $f(z) = h(z) + \overline{g(z)},$ is called sense-preserving at $z_0$ if the Jacobian $J_f(z) > 0$ for every $z$ in some punctured neighborhood of $z_0.$ We also say that $f$ is sense-reversing if $\overline{f}$ is sense-preserving at $z_0.$ If  $f$ is  neither sense-preserving nor sense-reversing at $z_0,$  then $f$  is said to be a  singular polynomial at $z_0.$
  \end{definition}
 \begin{theorem}\cite{kirwan1992complex} $\label{II0}$
  Let $f$ and $g$ be relatively prime polynomials in the real variables $x$ and $y$ with real coefficients, and let $\mathrm{deg}h = n$ and $\mathrm{deg}g = m.$ Then the two algebraic curves $f(x, y) = 0$ and $g(x, y) = 0$ have at most $mn$ points in common.
  \end{theorem}
   As an immediate consequence of Bezout's theorem in the plane, the harmonic polynomial $f(z)= h(z)+\overline{g(z)}$ has at most $n^2$  zeros where $\mathrm{deg}h(z)=n> m= \mathrm{deg}g.$ It is well known that if $f$ is a complex valued harmonic function that is locally univalent in a domain $\mathscr{D} \subset \mathbb{C},$ then its Jacobian, $J_f(z)= |h'|^2-|g'|^2$ never vanish for all $z \in \mathscr{D}$(see \cite{lewy1936non}).  As a result, a complex valued harmonic function $f(z)=h(z)+\overline{g(z)}$ is locally univalent and sense-preserving if and only if $h'(z) \neq 0$ and $|\omega (z)|<1,$ where $\omega (z)$ is a  dilatation  function of $f$ defined by $ \omega (z) = \frac{g'(z)}{h'(z)}.$ Note that dilatation is a measure of how a harmonic function is far from being analytic. For instance, the dilatation of analytic \bigskip function is zero. 
  \begin{theorem}\cite{wilmshurst1998valence}$\label{II1}$
If $f(z)= h(z)+\overline{g(z)} $ is a complex-valued harmonic polynomial such that $\mathrm{deg}h =n>m=\mathrm{deg}g$ and $\lim_{z \rightarrow \infty}f(z) = \infty,$ then $f(z)$ has at most $n^2$ zeros.
\end{theorem}
 The upper bound follows from applying Bezout'a theorem and the lower bound is based on the generalized argument principle and is sharp for each $m$ and $n.$ The upper bound is sharp which was shown by Wilmshurst \cite{wilmshurst1998valence} and this upper bound is sharp in general. For instance,  $\left(\frac{1}{i}\right)^n Q\left(iz+\frac{1}{2} \right) $ is a polynomial with $n^2$  zeros where $Q(z)=z^n+(z-1)^n+i\overline{z}^n -i(\overline{z}-1)^n.$ But it is natural to ask whether or not it can  \bigskip be improved for some interesting special classes of polynomials.\\
\section{Curves formed by vanishing discriminant}
\begin{figure}[!h]%
    \centering   
  \caption{Discriminantal Curve of harmonic polynomial $\overline{z}^4=z^2+az+b$}
  \label{Discriminantal} 
 \includegraphics[width=8cm]{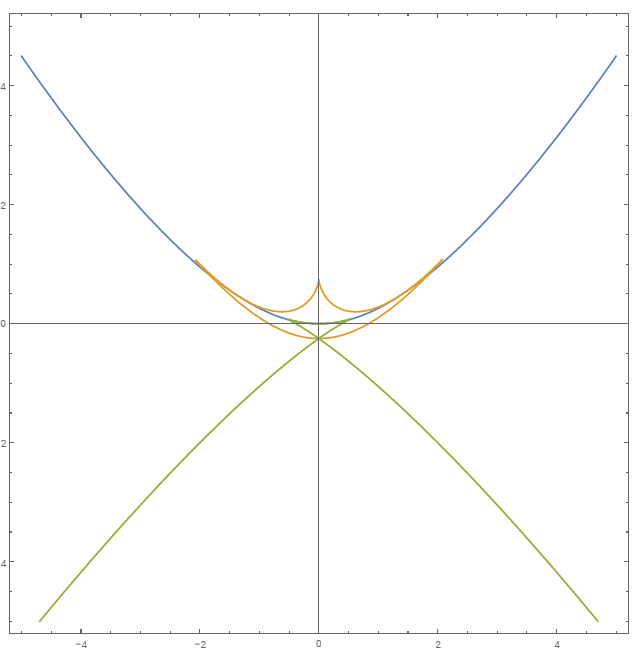}%
\end{figure}

$\label{c}$
Under this section, we show that the vanishing discriminant produces  four parabolic curves and one swallowtail catastrophe. Also by determining a concrete polynomial for each bounded and unbounded connected components of the intersecting curves, we show that the number of zeros of the family of harmonic polynomials of the type \bigskip $p(z)=\overline{z}^4-z^2-az-b$ cannot be more than $10.$\\
We are interested in counting the zeros of the harmonic polynomial equation 
\begin{equation} \label{1}
\overline{z}^4=z^2+az+b
\end{equation} for $a,b \in \mathbb{R}$ and also special types of curve are defined here. Now put $z=x+iy.$ Then  equation \ref{1} can be reduced to 
\begin{equation} \label{2}
 (x^4 + y^4 - 6x^2y^2 - x^2 + y^2 - ax - b) + i(4xy^2 - 4x^3 - 2x - a)=0
\end{equation}
Put $p_1(x,y):=x^4 + y^4 - 6x^2y^2 - x^2 + y^2 - ax - b$ and $p_2(x,y):=4xy^2 - 4x^3 - 2x - a.$ Equating both to zero, we have 
\begin{equation} \label{3}
x^4 + y^4 - 6x^2y^2 - x^2 + y^2 - ax - b=0
\end{equation} and 
\begin{equation} \label{4}
4xy^2 - 4x^3 - 2x - a=0.
\end{equation}
The Jacobian of $(p_1,p_2)$ denoted by $J_1$ given as follows.
 $$J_1=\frac{\partial p_1}{\partial x}\frac{\partial p_2 }{\partial y}-\frac{\partial p_1}{\partial y}\frac{\partial p_2}{\partial x}~~~~~~~~~~~~~~~~~~~~~~~~~~~~~~~~~~~~~~~~~~~~~~~~~~~~~~~~~~~~~~~~~~~~~~~~$$ \begin{equation} \label{5} 
 = 8xy(-a- 2x+4x^3-12xy^2)-(-2- 12x^2+4y^2)(2y-12x^2y+4y^3)
 \end{equation}
 By eliminating equations \ref{3}, \ref{4} and \ref{5} simultaneously, we get the discriminant of equation \ref{1}. Here, also we can find the discriminant for $J_1$ by singular program(software) and is given as follows.
 $$\label{6}
 \bigtriangleup _{a,b} = (36 a^2 + 125 a^4 - 144 b - 560 a^2 b + 384 b^2 - 256 b^3) (-6912 - 14976 a^2 + 13552 a^4 - 32616 a^6 +$$ $$91125 a^8 + 13824 b + 8832 a^2 b - 36000 a^4 b - 699840 a^6 b + 66816 b^2 + 566272 a^2 b^2 + 1645056 a^4 b^2 -$$ $$ 278528 b^3 - 1144832 a^2 b^3 - 152064 a^4 b^3 + 401408 b^4 + 589824 a^2 b^4 - 262144 b^5 + 65536 b^6).$$
Note that this discriminant can be calculated by singular program according to the following order.\\
$~~~~~~~~~\blacktriangleright ring ~~r=0,~(a,b,c,x,y),~dp;\\
~~~~~~~~~\blacktriangleright poly~ p_1=x^4 + y^4 - 6*x^2*y^2 - x^2 + y^2 - a*x - b;\\
~~~~~~~~~\blacktriangleright poly~ p_2=4*x*y^2 - 4*x^3 - 2*x - a=0;\\
~~~~~~~~~\blacktriangleright poly~ J = diff(p_1,x)*diff(p_2,y) - diff(p_2,x)*diff(p_1,y);\\
~~~~~~~~~\blacktriangleright ideal ~i = p_1,~p_2,~J;\\
~~~~~~~~~\blacktriangleright ideal ~j = std(i);\\
~~~~~~~~~\blacktriangleright ideal ~k = eliminate(j,x*y);\\
~~~~~~~~~\blacktriangleright k;\\ $ 
\bigskip Here, the final value, $k$ is the desired discriminant value.\\
 Similarly, one can calculate the jacobian and discriminant for the other case, that is for $y=0.$ By taking 
 \begin{equation} \label{7}
p_1:=x^4 + y^4 - 6x^2y^2 - x^2 + y^2 - ax - b
\end{equation} and 
\begin{equation} \label{8}
p_2*:=y,
\end{equation} one can calculate the second jacobian $J_2$ of $(p_1 , p_2*)$ as, 
\begin{equation} \label{9}
J_2=\frac{\partial p_1}{\partial x}\frac{\partial p_2 }{\partial y}-\frac{\partial p_1}{\partial y}\frac{\partial p_2}{\partial x} = -a - 2 x + 4 x^3 - 12 x y^2.
\end{equation} and the corresponding second discriminant for $J_2$ as 
\begin{equation}
\bigtriangleup _{a,b}=-4 a^2 + 27 a^4 + 16 b - 144 a^2 b + 128 b^2 + 256 b^3.
\end{equation}
Next, let us have the following notations:
$$\bigtriangleup_1 = 36 a^2 + 125 a^4 - 144 b - 560 a^2 b + 384 b^2 - 256 b^3;~~~~~~~~~~~~~~~~~~~~~~~~~~~~~~~~~~~~~~~~~~~~~~~~~~$$
$$\bigtriangleup_2 = -6912 - 14976 a^2 + 13552 a^4 - 32616 a^6 + 91125 a^8 + 
   13824 b + 8832 a^2 b - 36000 a^4 b$$ $$ - 699840 a^6 b + 66816 b^2 + 
   566272 a^2 b^2 +  1645056 a^4 b^2 - 278528 b^3 - 1144832 a^2 b^3 - $$$$
   152064 a^4 b^3 +  401408 b^4 + 589824 a^2 b^4 - 262144 b^5 + 
   65536 b^6;~~~~~~~~~~~~~~~~~~~~~~~~~~$$ and
   $$\bigtriangleup_3= -4 a^2 + 27 a^4 + 16 b - 144 a^2 b + 128 b^2 + 256 b^3~~~~~~~~~~~~~~~~~~~~~~~~~~~~~~~~~~~~~~~~~~~~~~~~~~~~~~~~~~~~$$ 
   Note that a polynomial is with zeros of multiplicity at least two if and only if its discriminant vanishes. Therefore, we focus our attention to the contour plot of the factored \bigskip discriminant(function of coefficient). \\The contour plot, which is a graphical technique for representing a surface, of 
   \begin{equation}
   \bigtriangleup_1=\bigtriangleup_2=\bigtriangleup_3=0
   \end{equation} on the interval $[-5,5]$ is given below as in Figure \ref{Discriminantal}. Equivalently, we find the curves formed by the following simultaneous equations.
   \begin{equation}
   \begin{cases}
   36 a^2 + 125 a^4 - 144 b - 560 a^2 b + 384 b^2 - 256 b^3=0\\ \\
   -6912 - 14976 a^2 + 13552 a^4 - 32616 a^6 + 91125 a^8 +\\ 
   13824 b + 8832 a^2 b - 36000 a^4 b$$ $$ - 699840 a^6 b + 66816 b^2 + \\
   566272 a^2 b^2 +  1645056 a^4 b^2 - 278528 b^3 - 1144832 a^2 b^3 - \\
   152064 a^4 b^3 +  401408 b^4 + 589824 a^2 b^4 - 262144 b^5 + 
   65536 b^6=0\\ \\
   -4 a^2 + 27 a^4 + 16 b - 144 a^2 b + 128 b^2 + 256 b^3=0
   \end{cases}
\end{equation} Here, it is recommendable to minimize the intervals of parameters to see all connected components. In Mathematica, we insert the following to sketch the curve:
$ContourPlot[{\bigtriangleup_1==0,\bigtriangleup_2==0,\bigtriangleup_3==0},{a,-5,5},{b,-5,5},PlotPoints \mapsto 200,\\ Axes \mapsto True]$

\subsection{Interpretation of the Curves formed by a Vanishing Discriminant}
\begin{figure}[htp]

\centering
\includegraphics[width=.3\textwidth]{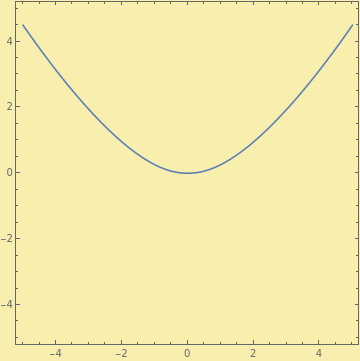}\hfill
\includegraphics[width=.3\textwidth]{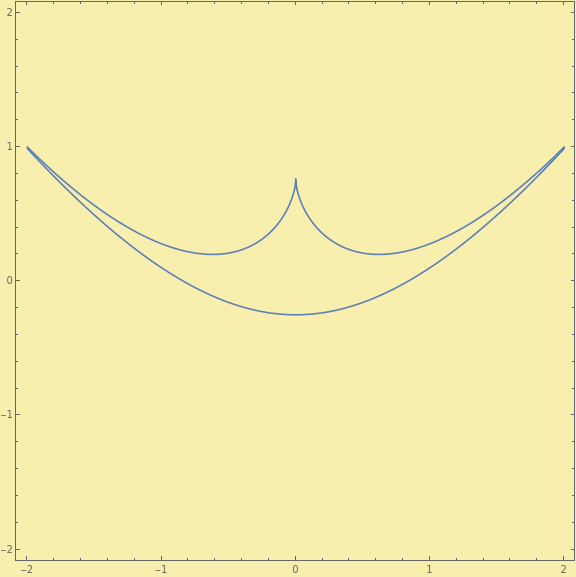}\hfill
\includegraphics[width=.3\textwidth]{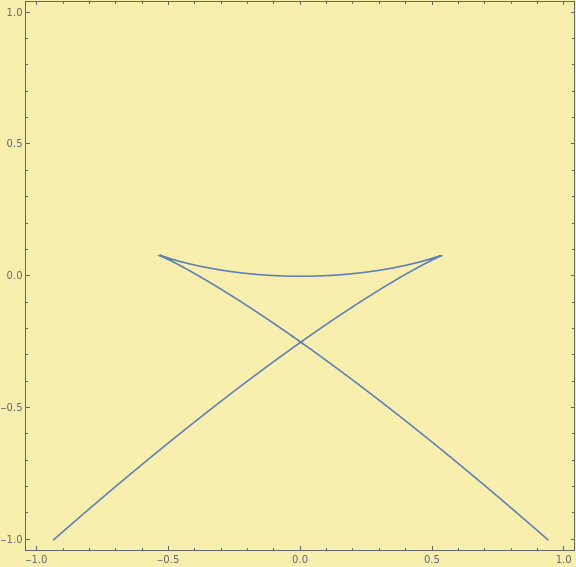}
\caption{Three curves(parabolic, teacup and swallowtail) formed by a vanishing discriminant}
\label{Three curves}

\end{figure}

Separately, we have the following three curves formed by a vanishing discriminant. As it can be seen from Figure \ref{Three curves} we have the following facts.
\begin{enumerate}
\item \textbf{\underline{The Left Corner:}} The curve at the left corner is a parabolic curve and is formed by $\bigtriangleup _1.$ 
\item \textbf{\underline{The Middle:}} The middle one is also a set of three parabolic curve coming together to form a special type of curve looks like the teacup and a bright curve(caustics) with its image in the teacup. Here if we assume that the lower big curve is a mirror, one is the image of the other and vice-versa. The center of a mirror(a curve of a teacup) is on the $b-axis$ and a bright curve is either to the left or to the right side of $b-axis.$ Actually this curve is formed due to $\bigtriangleup_2.$
\item \textbf{\underline{The Right Corner:}} The right corner curve is called a Swallowtail Catastrophe. It is formed due to $\bigtriangleup_3$.
\end{enumerate}
\section{Exact number of zeros in bounded and unbounded components}$\label{e}$
In this section, we put the lower and upper bound of the number of zeros of complex-valued harmonic polynomial families of the type $\overline{z}^4-(z^2+az+b)$ for  $a,b \in \mathbb{R}$ in each bounded and unbounded connected components formed by using a concrete \bigskip polynomial for each.\\
Note that the picture tells us that there are 12 connected components in the complement to the curve in the (a,b)-plane.  There are four unbounded components, namely $U_1, U_2, U_3$ and $U_4$ (Figure \ref{Unbounded components $U_1,...U_4$}), and eight bounded components, namely $B_1, B_2,...B_8$ (Figure \ref{Bounded components $B_1,...B_8$}). Now we have to pick one point in each of these components and we will get 12 concrete harmonic polynomials. We should then calculate their zeros in Mathematica and this will give a complete answer about the number of zeros of harmonic polynomials for this specific family. This is because this number stays the  same within each such connected component of the complement \bigskip to the curve. \\
Let us consider now the first unbounded component named $\mathrm{U_1}.$ Similarly we can do for the other components.\\
\textbf{\underline{In $\mathrm{U_1}$:}}\\

 Pick $(4,0).$ Then we solve the simultaneous equation $$x^4 + y^4 - 6x^2y^2 - x^2 + y^2 - ax - b = y(4xy^2 - 4x^3 - 2x - a)=0$$ with $a=4$ and $b=0.$ \\ Then by using Mathematica we find $6$ roots in this unbounded  component which are \bigskip given by\\ $ 0+0i, 1. 79632+0i,  -1. 11564 -0. 921028i ,  -1. 11564, + 0. 921028i, 0. 629282 -1. 57642i,$ \bigskip and  $ 0. 629282+ 1. 57642i .$
   \begin{figure}[!h]%
    \centering
  \caption{Unbounded components $U_1,...U_4$}
  \label{Unbounded components $U_1,...U_4$}
   \includegraphics[width=8cm]{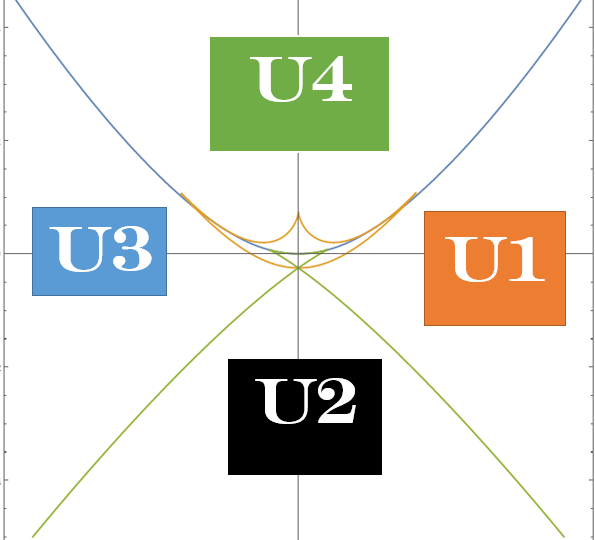}%
\end{figure}
\\
In similar fashion, by picking the points $(0,-4), (-4,0)$ and $(0,4)$ in $U_2, U_3$ and $U_4$ respectively we find $4,6$ and $4$ total number of roots in each component respectively. Hence, there are at most $6$ roots in unbounded components. Generally, the concrete harmonic polynomials and the corresponding number of zeros in each unbounded component are summarized as follow: 
\begin{table}[H]
\caption{The concrete polynomials in unbounded components}\label{UT}
\begin{center}
\begin{tabular}{|p{2cm}|p{5.5cm}|p{2cm}|p{2cm}|} \hline
Components & Concrete Polynomial &Size & Remarks \\ \hline
$U_1$ & $x^4 + y^4 - 6x^2y^2 - x^2 + y^2 - 4x - y(4xy^2 - 4x^3 - 2x - 4)i$ & 6 & distinct \\ \hline
$U_2$ & $x^4 + y^4 - 6x^2y^2 - x^2 + y^2 +4 - y(4xy^2 - 4x^3 - 2x )i$ & 4 & distinct \\ \hline
$U_3$ & $x^4 + y^4 - 6x^2y^2 - x^2 + y^2 + 4x - y(4xy^2 - 4x^3 - 2x +4)i$ & 6 & distinct \\ \hline
$U_4$ & $x^4 + y^4 - 6x^2y^2 - x^2 + y^2 - 4 - y(4xy^2 - 4x^3 - 2x )i$ &4& distinct \\ \hline
\end{tabular}
\end{center}
\end{table}
Therefore, the number of zeros of the family $f(z)=\overline{z}^4-z^2-az-b$ for $a,b \in \mathbb{R}$ is  bounded \bigskip below by $4$ and above by $10$ in all unbounded components.\\
Next, we consider the bounded components. As we did for unbounded components, we pick one point from each bounded component and find a concrete polynomial for  each \bigskip connected component. \\For instance, the point $(\frac{-1}{2},\frac{1}{15})$ is located in $B_5.$ The corresponding concrete polynomial for this bounded component is 
\begin{equation}
x^4 + y^4 - 6x^2y^2 - x^2 + y^2 +\frac{1}{2}x - \frac{1}{15} - y(4xy^2 - 4x^3 - 2x +\frac{1}{2})i.
\end{equation}  Using Mathematica, we find $10$ zeros in $B_5.$ 
\begin{figure}[!h]
    \centering
  \caption{Bounded components $B_1,...B_8$}
  \label{Bounded components $B_1,...B_8$}
   \includegraphics[width=10cm]{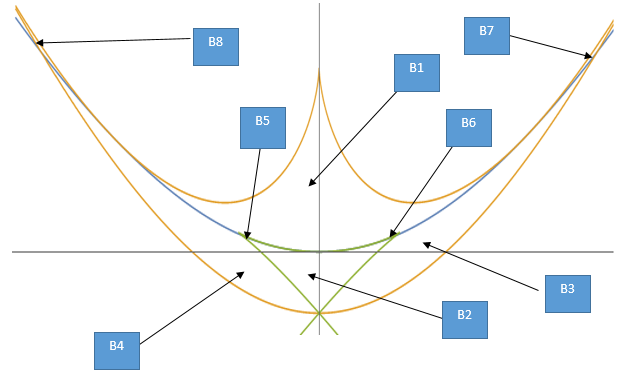}%
\end{figure}
Similarly, there are 8, 8, 6, 6, 10, 4 and 4 zeros in $B_1, B_2, B_3, B_4, B_6, B_7$ and $B_8$ respectively where each $B_i's$ are all bounded components and labeled as in the Figure \ref{Bounded components $B_1,...B_8$}.
\section{Conclusions} $\label{r}$
In this paper, we have seen that the curve formed on the $(a,b)$-plane are swallowtail catastrophe(due to $\bigtriangleup_3),$ a caustics and teacup curve(due to $\bigtriangleup_2$), swallowtail butterfly(due to $\bigtriangleup_1, \bigtriangleup_2$ and $\bigtriangleup_3$) and parabola(due to $\bigtriangleup_1).$ Generally, if we consider the curve formed by a vanishing discriminant, i.e, for $\bigtriangleup_1, \bigtriangleup_2$ and $\bigtriangleup_3$ all vanishing, the curve represents a catastrophe butterfly with  $\bigtriangleup_1$  is  the wing, $\bigtriangleup_3$ is is the leg, and $\bigtriangleup_2$ \bigskip is a beak and other body part.\\
The above analysis, in this paper, in both bounded and unbounded connected components gives a complete answer on the number of zeros of complex-valued harmonic polynomial of the form $p(z)=\overline{z}^4-z^2-az-b$ for $a,b \in \mathbb{R}.$ Thus, the number of zeros of $p(z)$ rises from $4$ to $10.$ An interesting example for the sharpness of this bound is a complex-valued harmonic polynomial with analytic part $-30z^2-15z-2$ and co-analytic part $30z^4.$ That means, a harmonic polynomial 
\begin{equation}\label{conclusion}
p(z)=-30z^2-15z-2+ 30\overline{z}^4
\end{equation} has $10$ number of roots. After expanding this equation, we use the following program to find all roots of $P.$ \\
$
~~~~~~~~~~p1[x_, y_] := x^4 + y^4 - 6*x^2*y^2 - x^2 + y^2 +\frac{1}{2}*x - \frac{1}{15};\\
~~~~~~~~~~p2[x_, y_] :=4*x*y^2 - 4*x^3 - 2*x +\frac{1}{2}\\
~~~~~~~~~~SOL = Solve[\{p1[x, y] == 0, y*p2[x, y] == 0\}, \{x, y\}, Reals]\\
~~~~~~~~~~N[SOL]\\
~~~~~~~~~~Length[\%] $\\
This program answers on the number of roots of equation $(\ref{conclusion})$ and the following \bigskip output is obtained showing that the total number of zeros is $10$ namely:\\
$-1.20812 + 0i, 0.254355+0i, 0.374583+0i, 0.579178+0i, -0.499333-0.999834i,-0.499333+0.999834i, 0.227692-0.0534519i, 0.227692+0.0534519i, 0.498134-0.705125i,$ and $ 0.498134+0.705125i$
\section*{For Further Investigation}
\begin{enumerate}
\item The authors' intention for the future work  is to consider the application of a complex-valued harmonic polynomials of the type $\overline{z}^4=z^2+az+b$ for $a,b \in \mathbb{R}$ in Bio-mathematics by combining all curves formed due to vanishing discriminant that forms a  butterfly catastrophe.
\item The curve formed due to $\bigtriangleup_2$ in this study is also of our interest to continue with. It is possible to extend this polynomial to caustics study by fixing a light source and bright curve. 
\end{enumerate}
\section*{Acknowledgments}
We would like to thank Stockholm University, Addis Ababa University, Simons Foundation,and ISP(International Science Program) at department level for providing us opportunities and financial support. 
\section*{Declaration of Interest of Statement}
The authors declare that there are no conflicts of interest regarding the publication of this paper.

\end{document}